\documentclass[5p,times,lefttitle,twocolumn]{elsarticle}

\def\figpath{figures}

\usepackage[T1]{fontenc}
\usepackage[utf8]{inputenc}
\usepackage{amsmath,amssymb,amsthm}
\usepackage{dsfont}
\usepackage{mathtools}
\usepackage{empheq}
\usepackage[binary-units]{siunitx}
\DeclareSIUnit\crate{C}
\usepackage{algorithm}
\usepackage{algpseudocode}
\usepackage{booktabs}
\usepackage{longtable}
\usepackage{doi}
\usepackage[capitalize, noabbrev]{cleveref}
\crefname{subsection}{Subsection}{Subsections}
\crefname{equation}{Eq.}{Eq.}

\usepackage{hyperref}
\hypersetup{
  hidelinks=true
}

\usepackage{multicol}

\let\oldbibitem\bibitem
\def\bibitem{\vfill\oldbibitem}

\journal{Examples and Counterexamples}

\biboptions{sort&compress}

\usepackage{amsmath,amssymb}%
\usepackage{dsfont}%
\usepackage{amsfonts}%
\usepackage{amsthm}%

\makeatletter
\newcommand{\@defs@vec}[1]{\boldsymbol{#1}}
\newcommand{\@defs@tens}[1]{\mathbf{#1}}

\DeclareMathOperator{\RelTol}{RelTol}
\DeclareMathOperator{\AbsTol}{AbsTol}
\DeclareMathOperator{\tr}{tr \!}

\newcommand{\Norm}[3]{\big|\!\big| #1 \big|\!\big|_{#2}^{#3}}

\newcommand{\de}{\,\mathrm{d}}
\newcommand{\ptl}{\partial}

\newcommand{\IR}{\mathbb{R}}
\newcommand{\IRp}{\mathbb{R}_{\geq0}}

\newcommand{\Ctensor}{\mathds{C}}

\newcommand{\vect}[1]{\@defs@vec{#1}{}}
\newcommand{\tens}[1]{\@defs@tens{#1}{}}

\newcommand{\trp}{\textsf{T}}

\newcommand{\grad}{\boldsymbol{\nabla}}

\newcommand{\tfinal}{t_\text{end}}

\newcommand{\OmegaL}{\Omega_0}
\newcommand{\OmegaLP}{\Omega_{0,\text{P}}}
\newcommand{\OmegaLS}{\Omega_{0,\text{S}}}
\newcommand{\PP}{\text{P}}
\newcommand{\SEI}{\text{S}}

\newcommand{\lagg}{\text{lag}}
\newcommand{\logg}{\text{log}}

\newcommand{\npo}{{n+1}}
\newcommand{\halb}{\frac{1}{2}}

\newcommand{\dualreduction}{\!:\!}

\newcommand{\tin}{\text{in }}
\newcommand{\ton}{\text{on }}

\newcommand{\ch}{\text{ch}}

\newcommand{\el}{\text{el}}
\newcommand{\cmax}{c_{\text{max}}}

\newcommand{\pl}{\text{pl}}
\newcommand{\gradL}{\grad_0}
\newcommand{\divgL}{\grad_0 \!\cdot\!}
\newcommand{\something}{\boldsymbol{\cdot}}
\newcommand{\minus}{\scalebox{0.75}[1.0]{$-$}}
\newcommand{\minusone}{{\minus 1}}
\newcommand{\normal}{\vect{n}}

\newcommand{\ext}{\text{ext}}
\newcommand{\devi}{\text{dev}}

\newcommand{\q}{\quad}

\newcommand{\fdg}{\,:\,}

\newcommand{\conb}{\OB{c}}

\newcommand{\OB}[1]{\mkern 1.5mu\overline{\mkern-1.5mu#1\mkern-1.5mu}\mkern
  1.5mu}

\makeatother

\newtheoremstyle{boldremark}
{\dimexpr\topsep/2\relax} 
{\dimexpr\topsep/2\relax} 
{}          
{}          
{\bfseries} 
{.}         
{.5em}      
{}          

\theoremstyle{boldremark}

\begin{document}

  \begin{frontmatter}

    \title{%
      Comparison of Different Elastic Strain Definitions
      for
      Largely
      Deformed SEI
      of Chemo-Mechanically Coupled Silicon Battery Particles
    }

    \author[ianm]{R. Schoof\corref{correspondingauthor}}
    \cortext[correspondingauthor]{Corresponding author}
    \ead{raphael.schoof@kit.edu}

    \author[tvt]{G. F. Castelli}

    \author[ianm]{W. D\"orfler}

    \address[ianm]{Karlsruhe Institute of Technology~(KIT),
      Institute for Applied and Numerical Mathematics,
      Englerstr.~2, 76131~Karlsruhe, Germany}

    \address[tvt]{%
      Fraunhofer~Institute~for~Industrial~Mathematics~ITWM,
      Fraunhofer-Platz~1,
      67663~Kaiserslautern,~Germany
    }



\begin{abstract}
  Amorphous silicon is a highly promising anode material for next-generation
  lithium-ion batteries.
  Large volume changes of the silicon particle
  have
  a critical effect on the surrounding
  solid-electrolyte interphase (SEI) due to repeated fracture and healing
  during
  cycling.
  Based on a thermodynamically consistent chemo-elasto-plastic continuum model
  we investigate the stress development inside the particle and the SEI.
  Using the example of a particle with SEI,
  we apply a higher order finite element method
  together with
  a variable-step, variable-order time integration scheme
  on a nonlinear system of partial differential equations.
  Starting from a single silicon particle setting,
  the surrounding SEI is added in a first step
  with the
  typically used elastic Green--St-Venant (GSV) strain
  definition for a purely elastic deformation.
  For this type of deformation,
  the definition of the elastic strain is crucial
  to get reasonable
  simulation results.
  In case of the elastic GSV strain,
  the simulation aborts.
  We overcome the simulation failure by using the definition of the logarithmic
  Hencky strain.
  However, the particle remains unaffected by the elastic strain definitions
  in the particle domain.
  Compared to GSV, plastic deformation with the Hencky strain is straightforward
  to take into account.
  For the plastic SEI deformation,
  a rate-independent and a rate-dependent plastic deformation are newly
  introduced
  and numerically
  compared for
  three half cycles for the example of a radial symmetric particle.
\end{abstract}


    \begin{keyword}
      finite elements \sep
      numerical simulation \sep
      solid-electrolyte interphase (SEI) \sep
      finite elastic strain \sep
      (visco-)plasticity
      \MSC[2020]
      74C15 \sep
      74C20 \sep
      74S05 \sep
      65M22 \sep
      90C33
    \end{keyword}

  \end{frontmatter}



\section{Introduction}
\label{sec:introduction}

Nowadays, the advancement of lithium-ion batteries is one key research in terms
of climate
change~\cite{sarmah2023recent}.
Alongside all solid-state batteries
as promising candidates for future battery
types~\cite{landstorfer2011advanced},
especially, amorphous silicon (aSi) is superior in terms of the nearly tenfold
theoretical energy density compared to state of the art graphite
anodes
for classical battery types~\cite{zhao2019review}.
Unfortunately, this capacity increase is accompanied by a volume change up to
300\%
crucially effecting the surrounding solid-electrolyte interphase (SEI) during
swelling
and
shrinking~\cite{zhao2019review, kolzenberg2022chemo-mechanical}.
The SEI features elastic and plastic deformation before it can break up and
heal again during repeated cycling~\cite{kolzenberg2022chemo-mechanical}.
In addition,
the voltage hysteresis of silicon is still an ongoing research topic
because it influences battery lifetime and performance
and makes further mechanical investigation of the SEI
unavoidable~\cite{kobbing2024voltage}.
Recent measurements of silicon anodes show large overpotentials with slow
relaxation
which could confirm the mechanical stress influence of the SEI and could
explain the voltage hysteresis~\cite{wycisk2023challenges}.

In this paper,
we apply modern numerical techniques on the extended underlying
system of
equations, based on an efficient adaptive algorithm,
compare~\cite{castelli2021efficient}
and~\cite[Sect.~4.1]{castelli2021numerical}.
Furthermore,
we discuss the influence of the definition of the
elastic strain tensor for the SEI domain
on the numerical method,
starting with a single particle
setting~\cite{castelli2021efficient}.
Due to the large volume change, a finite deformation approach is needed for the
particle as well as for the SEI domain.
Firstly, a purely elastic case for the SEI is considered with two different
definitions for the elastic strain tensor. We use the definition of the
(Lagrangian)
Green--St-Venant (GSV)
strain tensor and the logarithmic (Hencky) strain tensor.
In the situation of a single particle, both definitions lead to
quantitative
similar
numerical results~\cite{schoof2023efficient, zhang2018sodium}.
The strain definition plays a significant role for the SEI domain whether the
numerical simulation will abort or not.
To the best of the authors' knowledge, this has not yet been documented in the
literature.
Using the logarithmic strain approach, we
can comfortably add a
rate-independent and a rate-dependent plastic deformation for the SEI without
increasing the system of equation for the plastic part of the deformation
gradient~\cite{kolzenberg2022chemo-mechanical,
schoof2023efficient}.
The rate-dependent case results is a typical stress-overrelaxation, which is
also
found in measurements of the electric field~\cite{wycisk2023challenges}
and can confirm the hypothesis of~\cite{kobbing2024voltage}.

The rest of the article is structured as follows:
in the next section,
we present and recap our used model approach and present
the different definitions for the elastic strain approaches.
In~\cref{sec:numerics}, we formulate our model equations and show important
steps towards a numerical solution.
\cref{sec:results} shows the failure of one of the used definitions of the
elastic strains and successful numerical simulations for the other
definition.
In addition, the latter one is extended with plastic and viscoplastic
deformation.
In the end, we sum up our main findings and close with an outlook.




\section{Theory}
\label{sec:theory}
In this section, we recall and summarize our used model equations for the
chemo-mechanically coupled
\textit{particle-SEI}
approach
based on the thermodynamically consistent theory
by~\cite{schoof2023efficient, kolzenberg2022chemo-mechanical,
castelli2021efficient}.
For the electrode particle,
we use a pure elastic deformation model with a Green--St-Venant (GSV) strain
tensor~\cite{kolzenberg2022chemo-mechanical}.
For the SEI, we introduce four different
approaches
for the deformation:
1.) a pure elastic with the GSV strain,
2.) a pure elastic (Lagrangian) logarithmic Hencky strain
as well as
3.) a plastic and
4.) a viscoplastic deformation approach,
whereby last two use the
logarithmic
Hencky
strain~\cite{schoof2023efficient, neff2016exponentiated, neff2016geometry}.

\textbf{Finite Deformation.}
Let~$\OmegaL \subset \IR^3$
be a bounded domain which represents the particle-SEI domain in the reference
(Lagrangian) configuration,
divided into a particle
domain~$\OmegaLP$
and a SEI
domain~$\OmegaLS$.
We denote the current (Eulerian) domain configuration without the
subscript~$0$.
With the
displacement~$\vect{u}$
and the
deformation~$\Phi\colon \IRp \times \OmegaL \rightarrow \Omega$,
$\Phi(t, \vect{X_0}) \coloneqq \vect{x} = \vect{X}_0 + \vect{u}(t, \vect{X}_0)$
from the reference configuration to the current
configuration~\cite[Sect.~2]{holzapfel2010nonlinear},
\cite[Sect.~8.1]{lubliner2006plasticity},
\cite[Ch.~VI]{braess2007finite},
\cite{castelli2021efficient, kolzenberg2022chemo-mechanical,
schoof2023efficient, schoof2023simulation, di-leo2015diffusion-deformation},
we can derive the deformation
gradient~$\tens{F} (\gradL \vect{u})\in \IR^{3,3}$
as~$\tens{F}
\coloneqq \ptl \Phi / \ptl \vect{X}_0
= \tens{Id} + \gradL \vect{u}$
with the identity tensor~$\tens{Id}$.
Following~\cite[Sect.~10.4]{bertram2021elasticity},
\cite[Sect.~8.2.2]{lubliner2006plasticity}
and~\cite{castelli2021efficient, di-leo2015diffusion-deformation},
we multiplicatively split the deformation
gradient~$\tens{F} = \tens{F}_\ch\tens{F}_\el\tens{F}_\pl$
into three parts: chemical, elastic and plastic part, respectively.
Note
that we have for each
domain~$\OmegaLP$
and~$\OmegaLS$
own deformation gradients~$\tens{F}_\PP$ and $\tens{F}_\SEI$.
We
consider only chemical and elastic deformation for the particle
domain~($\tens{F}_{\pl, \PP} = \tens{Id}$),
whereas only elastic and plastic deformations are considered for the SEI domain
(no calculation of any chemical species, $\tens{F}_{\ch, \SEI} = \tens{Id}$).
Furthermore, we identify two displacements:
the particle
displacement~$\vect{u}_\PP
\colon
\OB{\Omega}_{0, \PP, \tfinal} \rightarrow \IR^3$
and the SEI
displacement~$\vect{u}_\SEI
\colon
\OB{\Omega}_{0, \SEI, \tfinal} \rightarrow \IR^3$
with the final simulation
time~$\tfinal > 0$
and~$\OB{\Omega}_{0, \square, \tfinal} \coloneqq
[0, \tfinal] \times \OB{\Omega}_{0, \square}$,
$\square \in \{ \PP, \SEI\}$.
The chemical part arises due to lithium insertion in the particle domain and is
stated
as~$\tens{F}_{\ch, \PP} (\conb_\PP(t, \vect{X}_{0, \PP}))
= \tens{F}_\ch(\conb)
= \lambda_\ch (\conb) \tens{Id}
= \sqrt[3]{1+ v_\text{pmv} \cmax \conb} \tens{Id}$
with the partial molar
volume~$v_\text{pmv}> 0$
of lithium inside
aSi
and the normalized lithium concentration~$\conb = c / \cmax \in [0, 1]$
of the lithium concentration~$c
\colon$
$\OB{\Omega}_{0, \PP, \tfinal}$ $\rightarrow$ $[0, \cmax]$
with respect to the maximal concentration of
aSi,
$\cmax$,
in the
reference configuration.
Using the chemical
part~$\tens{F}_\ch$
and the plastic
part~$\tens{F}_\pl(t, \vect{X}_0)$,
applied as internal
variable~\cite{schoof2023efficient},
the elastic part can be computed
as~$\tens{F}_{\el, \PP} (\conb, \gradL \vect{u}_\PP)
=
\tens{F}_{\el} (\conb, \gradL \vect{u}_\PP)
= \lambda_\ch^{\minusone} \tens{F}$
and
$\tens{F}_{\el, \SEI} (\gradL \vect{u}_\SEI)
=
\tens{F}_{\el} (\gradL \vect{u}_\SEI)
= \tens{F} \tens{F}_\pl^{\minusone}$
with the
gradient~$\gradL$
in the respective domain.
If it is clear from the context which domain is considered,
the index~$\PP$
or~$\SEI$
is omitted for reasons of better readability.

\textbf{Free Energy.}
Following the thermodynamically consistent material
approach~\cite{kolzenberg2022chemo-mechanical, schoof2023efficient},
we additively decompose the Helmholtz free
energy~$\psi$
into a chemical
part~$\psi_\ch$
and a mechanical
part~$\psi_\el$ at constant temperature in the reference configuration:
\begin{align}
  \psi(\conb, \gradL \vect{u}, \tens{F}_\pl)
  =
  \psi_\ch(\conb)
  +
  \psi_\el(\conb, \gradL \vect{u}, \tens{F}_\pl).
\end{align}
For the respective domains,
we
have~$\psi_\PP(\conb, \gradL \vect{u}_\PP)
=
\psi_\ch(\conb)
+
\psi_{\el, \PP}(\conb, \gradL \vect{u}_\PP)$
and~$\psi_\SEI(\gradL \vect{u}_\SEI, \tens{F}_\pl)
=
\psi_{\el, \SEI}(\gradL \vect{u}_\SEI, \tens{F}_\pl)$.
The chemical free energy density is defined by an experimental open-circuit
voltage~(OCV)
curve~$U_\text{OCV}$~\cite{kolzenberg2022chemo-mechanical, schoof2023efficient,
schoof2023simulation}
\begin{align}
  \rho_0\psi_\ch(\conb)
  =
  \minus
  \cmax
  \int_0^{\conb} \mathrm{Fa} \, U_\text{OCV}(z) \de z
\end{align}
with the mass
density~$\rho_0$
of aSi in the reference configuration
and the Faraday
constant~$\mathrm{Fa}$.
The mechanical free energy density is stated
for both the particle and the SEI domain
by a linear elastic
approach~\cite{kolzenberg2022chemo-mechanical, castelli2021efficient,
schoof2023efficient}
\begin{align}
  \rho_0 \psi_\el (\conb, \gradL \vect{u}, \tens{F}_\pl)
  =
  \halb \tens{E}_\el(\conb, \gradL \vect{u}, \tens{F}_\pl) \dualreduction
  \Ctensor[\tens{E}_\el]
\end{align}
with the elastic strain tensor~$\tens{E}_\el$,
the constant, isotropic stiffness tensor~$\Ctensor$ of aSi
and
$\Ctensor[\tens{E}_\el]
= \lambda \tr \left(\tens{E}_\el\right) \tens{Id}
+
2G \tens{E}_\el$.
The first and second Lam\'e
constants~$\lambda = 2G \nu / (1-2\nu)$
and~$G = E / \big( 2(1+\nu)\big)$
depend on the Young's modulus~$E$
and the Poisson's ratio~$\nu$.
The parameters for the particle and the SEI material are specified
in~\cref{sec:results}.

Our main concern in this paper is the definition of the elastic strain tensor.
In~\cite{schoof2023efficient},
two different definitions for
the elastic strain tensor are stated:
the GSV strain tensor, also called \textit{the} Lagrangian strain
tensor~\cite[Sect.~8.1]{lubliner2006plasticity}
\begin{align}
  \tens{E}_{\el, \lagg}
  = \tens{E}_{\el, \text{GSV}}
  = \halb \left(\tens{C}_\el - \tens{Id}\right)
  = \halb \left(\tens{F}_\el^\trp \tens{F}_\el - \tens{Id}\right),
\end{align}
and
the (Lagrangian) logarithmic Hencky strain tensor
\begin{align}
  \tens{E}_{\el, \logg}
  = \ln \left(\tens{U}_\el\right)
  = \ln \left(\sqrt{\tens{C}_\el} \right)
  = \sum_{\alpha = 1}^{3} \ln \left(\sqrt{\eta_{\el, \alpha}}\right)
    \vect{r}_{\el, \alpha} \otimes \vect{r}_{\el, \alpha}
\end{align}
with the
eigenvalues~$\eta_{\el, \alpha}$
and
eigenvectors~$\vect{r}_{\el, \alpha}$
of~$\tens{U}_\el$
being the unique, positive definite and symmetric right stretch part of the
unique polar decomposition
of~$\tens{F}_\el =
\tens{R}_\el\tens{U}_\el$~\cite[Sect.~2.6]{holzapfel2010nonlinear}.
The right Cauchy--Green
tensor~$\tens{C} = \tens{F}^\trp \tens{F}$
is also symmetric and positive
definite~\cite[Sect.~8.1]{lubliner2006plasticity}.
For the situation of considering an elastic particle without SEI only,
both elastic strain approach results have similar
outcomes~\cite{schoof2023efficient, zhang2018sodium}.
In this paper, we use for the particle domain the GSV strain.
For the SEI domain,
we will apply both strain approaches,
since especially the logarithmic strain is favorable to incorporate plastic
deformation~\cite{schoof2023efficient}.

\textbf{Chemistry.}
A generalized diffusivity
equation~\cite{anand2012cahn-hilliard-type, di-leo2015chemo-mechanics,
kolzenberg2022chemo-mechanical}
is used to describe the change of lithium concentration inside the reference
particle domain~$\OmegaLP$
\begin{align}
  \ptl_t c
  =
  \minus \divgL \vect{N}
\end{align}
with the lithium
flux~$\vect{N}(\conb, \gradL \vect{u}_\PP, \gradL \mu)
= \minus m (\conb, \gradL \vect{u}_\PP)
\gradL \mu
= \minus D \left(\ptl_c \mu \right)^{\minusone} \gradL \mu$,
the scalar
mobility~$m > 0$ for the applied isotropic case,
the diffusion
coefficient~$D > 0$ for lithium atoms in aSi.
The chemical
potential~$\mu
\colon \OB{\Omega}_{0, \PP,\tfinal} \rightarrow \IR$
is given as the partial derivative of the free energy density with respect
to~$c$ as $\mu = \ptl_c (\rho_0 \psi)$~\cite{schoof2023simulation,
schoof2024residual}:
\begin{align}
  \mu_{\lagg}
  &
  =
  \minus \mathrm{Fa} \, U_\text{OCV}
  -
  \frac{v_\text{pmv}}{3 \lambda_\ch^5} \tens{F}^\trp \tens{F} \dualreduction
  \Ctensor[\tens{E}_\el],
  \\
  \mu_{\logg}
  &
  =
  \minus \mathrm{Fa} \, U_\text{OCV}
  -
  \frac{v_\text{pmv}}{3 \lambda_\ch^3} \tr \left(\Ctensor[\tens{E}_\el]
     \right)
  .
\end{align}

At the surface of the particle domain, we apply a uniform and constant
external lithium flux~$N_\ext$.
The sign of this flux is positive or negative,
depending on lithium insertion and extraction, respectively,
and is expressed in terms of the charging rate
(\si{C}-rate).
The state of charge (SOC) links the simulation time, the external lithium flux
and the initial concentration
via~$\text{SOC}(t) = \conb_0 + N_\ext t$.
For more information about $N_\ext$, the \si{C}-rate and the SOC,
we refer
to~\cite{schoof2023efficient, kolzenberg2022chemo-mechanical,
castelli2021efficient}
and the references cited therein.

\textbf{Elastic and Inelastic Deformation.}
In both domains, the static balance of linear of momentum is used to consider
the
deformation~\cite{kolzenberg2022chemo-mechanical, castelli2021efficient,
schoof2023efficient, schoof2023simulation}:
\begin{align}
  \vect{0}
  =
  \divgL \tens{P}_\PP(\conb, \gradL \vect{u}_\PP),
  \qquad
  \vect{0}
  =
  \divgL \tens{P}_\SEI(\gradL \vect{u}_\SEI, \tens{F}_\pl)
\end{align}
with the first Piola--Kirchhoff
tensor~$\tens{P} \in \IR^{3,3}$,
thermodynamically consistent derived
as~$\tens{P} = \ptl_\tens{F}(\rho_0 \psi)$.
This results
in~$\tens{P}_{\PP, \lagg}
= \lambda_\ch^{-2} \tens{F} \Ctensor_\PP[\tens{E}_{\el, \lagg}]$
for the particle domain
and
in~$\tens{P}_{\SEI, \lagg}
= \tens{F}
\left(\tens{F}_\pl^{\minusone} \right)^{\trp}\tens{F}_\pl^{\minusone}
\Ctensor[\tens{E}_{\el, \lagg}]$
or
$\tens{P}_{\SEI, \logg}
= \tens{F}
\left(\tens{F}_\el^\trp \tens{F}_\el^{\vphantom{\trp}}\right)^{\minusone}
\left(\tens{F}_\pl^{\minusone} \right)^{\trp}\tens{F}_\pl^{\minusone}
\Ctensor[\tens{E}_{\el, \logg}]$
for the SEI domain
depending on the definition of the GSV strain tensor or the logarithmic strain,
respectively~\cite{kolzenberg2022chemo-mechanical, schoof2023efficient,
schoof2023simulation, kobbing2024voltage}.
The Cauchy stress $\boldsymbol{\sigma} \in \IR_{\text{sym}}^{3,3}$ in the
current configuration is defined
as~$\boldsymbol{\sigma}
= \tens{P} \tens{F}^\trp /
\det \left(\tens{F}\right)$~\cite[Sect.~3.1]{holzapfel2010nonlinear}
with~$\det \left(\tens{F}\right) >
0$~\cite[Sect.~2.4]{holzapfel2010nonlinear}.
For the inelastic deformation approach,
we base on the theory
in~\cite{schoof2023efficient}.
There are a rate-independent plastic approach with isotropic hardening
and a rate-dependent plastic approach developed and compared.
The special feature is the usage of a projector formulation mapping the
stresses onto the set of admissible
stresses~\cite{frohne2016efficient},
a concept also known as static
condensation~\cite{wilson1974static,
di-pietro2015hybrid}.
Following~\cite{schoof2023efficient},
we define for the SEI
domain~$\OmegaLS$
the Mandel stress~$\tens{M}_\SEI (\gradL \vect{u}_\SEI, \tens{F}_\pl)
= \tens{M} =\ptl_{\tens{E}_\el} (\rho_{0, \SEI} \psi_{\el, \SEI})
= \Ctensor_\SEI[\tens{E}_{\el}]$.
For the rate-independent ideal plasticity ($\gamma^\text{iso} = 0$
in~\cite{schoof2023efficient, frohne2016efficient}),
we express the classical loading and unloading conditions via the
Karush--Kuhn--Tucker~(KKT)~\cite[Sect.~1.2.1]{simo1998computational},
\cite[Sect.~3.2]{lubliner2006plasticity},
\cite{kolzenberg2022chemo-mechanical}
as
\begin{align}
  F_\text{Y} \leq 0,
  \qquad
  \dot{\varepsilon}_\pl^\text{eq} \geq 0,
  \qquad
  F_\text{Y} \dot{\varepsilon}_\pl^\text{eq} = 0
\end{align}
with the yield
function~$F_\text{Y} (\gradL \vect{u}_\SEI, \tens{F}_\pl,
\varepsilon_\pl^\text{eq})
= \|\tens{M}^\text{dev}\| -
\sigma_\text{Y}$,
the deviatoric Mandel
stress~$\tens{M}^\text{dev}
= \tens{M} - 1/3 \tr\left( \tens{M} \right) \tens{Id}$,
the yield stress~$\sigma_\text{Y}$
and the accumulated equivalent inelastic
strain~${\varepsilon}_\pl^\text{eq}(t, \vect{X}_{0, \SEI}) \geq 0$.
To be consistent with the one dimensional tensile test, we rescale the yield
stress with the
factor~$\sqrt{2/3}$~\cite[Sect.~2.3.1]{simo1998computational}.
In the viscoplastic approach,
an evolution equation of the equivalent plastic strain, given by
\begin{subequations}
  \label{eq:viscoplatic_epsplver}
  \begin{empheq}[
    left={
      \dot{\varepsilon}_\pl^\text{eq} =
      \empheqlbrace}
    ]{alignat=3}
    &0,\label{eq:viscoplatic_epsplver_a}
    && \Norm{\tens{M}^{\devi}}{}{} \leq \sigma_\text{Y},
    \\
    &
    \dot{\varepsilon}_0
    \Bigg(
    \frac{\Norm{\tens{M}^{\devi}}{}{} - \sigma_\text{Y}}
    {\sigma_{\text{Y}^{*}}}\Bigg)^\beta,
    \q
    &&
    \label{eq:viscoplatic_epsplver_b}
    \Norm{\tens{M}^{\devi}}{}{} > \sigma_\text{Y},
  \end{empheq}
\end{subequations}
describes the plastic deformation
instead of the KKT
conditions~\cite[Sect.~1.7]{simo1998computational},
\cite{di-leo2015diffusion-deformation, schoof2023efficient}.
The positive-valued stress-dimensioned
constant~$\sigma_{\text{Y}^{*}}$,
the reference tensile
stress~$\dot{\varepsilon}_0$
and the measure of the strain rate sensitivity of the
material~$\beta$
are defined
in~\cref{sec:results}.




\section{Numerical Approach}
\label{sec:numerics}

\subsection{Problem Formulation}
\label{subsec:problem_statement}

Using Table~1
in~\cite{schoof2023efficient},
we arrive at the dimensionless initial boundary value problem with inequality
conditions:
find the normalized concentration~$c
\colon \OB{\Omega}_{0, \PP, \tfinal} \rightarrow
[0,1]$,
the chemical potential~$\mu
\colon \OB{\Omega}_{0, \PP, \tfinal} \rightarrow \IR$,
the displacements~$\vect{u}_\PP
\colon \OB{\Omega}_{0, \PP, \tfinal} \rightarrow \IR^3$
and~$\vect{u}_\SEI
\colon \OB{\Omega}_{0, \SEI, \tfinal} \rightarrow \IR^3$
satisfying
\vspace*{-0.12cm}
\begin{subequations}
  \label{eq:ibvp}
  \begin{empheq}[left=\empheqlbrace]{alignat=4}
    \partial_t c
    &= \minus \divgL \vect{N}
    (c, \gradL \vect{u}_\PP)
    && \qquad
    && \tin \Omega_{0, \PP, \tfinal},
    \label{eq:ibvp_a}
    \\
    \mu
    &= \partial_c (\rho_0
    \psi
    (c, \gradL \vect{u}_\PP))
    &&
    && \tin \Omega_{0, \PP, \tfinal},
    \label{eq:ibvp_b}
    \\
    \tens{0}
    &=
    \divgL \tens{P}_\PP
    (c, \gradL \vect{u}_\PP)
    &&
    && \tin \Omega_{0, \PP, \tfinal},
    \label{eq:ibvp_c}
    \\
    \tens{0}
    &=
    \divgL \tens{P}_\SEI
    (\gradL \vect{u}_\SEI, \tens{F}_\pl)
    &&
    && \tin \Omega_{0, \SEI, \tfinal},
    \label{eq:ibvp_d}
    \\
    F_\text{Y}
    \leq 0,\q&
    \dot{\varepsilon}_\pl^\text{eq} \geq 0, \q
    F_\text{Y}
    \dot{\varepsilon}_\pl^\text{eq} = 0
    &&
    && \tin \Omega_{0, \SEI, \tfinal},
    \label{eq:ibvp_e}
    \\
    \vect{N} \cdot \normal_\PP
    &= N_\ext
    &&
    && \ton \ptl \Omega_{0, \PP, \tfinal}
    \label{eq:ibvp_f}
    \\
    \vect{u}_\PP
    &= \vect{u}_\SEI
    &&
    && \ton \Gamma_{\text{inter}, \tfinal}
    \label{eq:ibvp_g}
    \\
    \tens{P}_\PP \cdot \normal_\PP
    &= \tens{P}_\SEI \cdot \normal_\PP
    &&
    && \ton \Gamma_{\text{inter}, \tfinal}
    \label{eq:ibvp_h}
    \\
    \tens{P}_\SEI \cdot \normal_\SEI
    &= \tens{0}
    &&
    && \ton \ptl \hat{\Omega}_{0, \SEI, \tfinal}
    \label{eq:ibvp_i}
    \\
    c(0, \something)
    &= c_0
    &&
    && \tin \Omega_{0, \PP},
    \label{eq:ibvp_j}
    \\
    \tens{F}_\text{pl}(0, \something)
    &= \tens{Id}
    &&
    && \tin \Omega_{0, \SEI},
    \label{eq:ibvp_k}
    \\
    \varepsilon_\pl^{\text{eq}}(0, \something)
    &= 0
    &&
    && \tin \Omega_{0, \SEI}
    \label{eq:ibvp_l}
  \end{empheq}
\end{subequations}
\vspace*{-0.04cm}
with the interface of particle and SEI domain~$\Gamma_{\text{inter}, \tfinal}
\coloneqq
[0, \tfinal] \times \Gamma_{\text{inter}}
=
[0, \tfinal] \times
(\ptl \Omega_{0, \PP} \cap \ptl \Omega_{0, \SEI})
$,
$\ptl\hat{\Omega}_{0, \SEI, \tfinal}
\coloneqq [0, \tfinal] \times (\ptl\Omega_{0, \SEI}
\textbackslash \Gamma_{\text{inter}})$,
the outward unit normal vector~$\normal_\PP = \minus \normal_\SEI$
and the external boundary~$\ptl \Omega_{0, \SEI, \tfinal}
\coloneqq
[0, \tfinal] \times \ptl \Omega_{0, \SEI}$.
We assume a constant lithium flux with changing sign due to cycling
at the
particle boundary as well as
equal displacements and stresses at the interface
boundary~$\Gamma_{\text{inter}, \tfinal}$
as well as no stresses at the SEI
boundary~$\ptl \Omega_{0, \SEI}$.
The initial concentration~$c_0$ is chosen to be constant and
rigid body motions are excluded with appropriate displacement boundary
conditions.
For the detailed treatment of the plastic part~$\tens{F}_\pl$
of the deformation gradient
and the equivalent plastic strain~$\varepsilon_\pl^\text{eq}$
as internal variables,
we refer to~\cite{schoof2023efficient}.
Finally, the desired quantities,
the Cauchy stresses~$\boldsymbol{\sigma}_\PP$ and $\boldsymbol{\sigma}_\PP$,
are computed in a postprocessing step
using the solution variables,
$\tens{F}$, $\tens{F}_\el$, $\tens{E}_\el$ and $\tens{P}$.

\subsection{Numerical Solution Procedure}
\label{subsec:num_sol_approach}

\textbf{Weak Formulation.}
Following~\cite{schoof2023efficient, castelli2021efficient,
frohne2016efficient},
we multiply with test functions,
integrate over the respective domain,
integrate by parts
and
formulate a weak primal mixed variational inequality for~\cref{eq:ibvp}.
Using a projector~$\tens{P}_\Pi$
onto the set of admissible
stresses~\cite{frohne2016efficient, suttmeier2010on},
we can reformulate the occurring saddle point problem as a primal formulation:
for given~$\tens{F}_\pl$
and~$\varepsilon_\pl^\text{eq}$
find solutions~$\lbrace c, \mu, \vect{u}_\PP, \vect{u}_\SEI \rbrace$
with~$c, \mu \in V$,
$\ptl_t c \in L^2(\Omega_{0}, \IR)$,
$\vect{u}_\PP \in \vect{V}_\PP^{*}$
and
$\vect{u}_\SEI \in \vect{V}_\SEI^{*}$
such that
\begin{subequations}
  \label{eq:weak_formulation_primal}
  \begin{empheq}[left=\empheqlbrace]{alignat=2}
    \big( \varphi,\partial_t c \big)
    &= \minus
    \big(
    \gradL \varphi,
    m
    (c, \grad_0 \vect{u}_\PP)
    \gradL \mu
    \big)
    - \big(\varphi, N_\ext \big)_{\Gamma_{\text{inter}}},
    \vspace{0.2cm}
    \label{eq:weak_formulation_primal_a}
    \\
    0 &= \minus
    \big(
    \varphi, \mu
    \big)
    + \big( \varphi, \ptl_c (
    \rho_0
    \psi_\ch(c))
    \nonumber
    \\
    &
    \quad \quad \quad \quad \quad
    + \ptl_c (
    \rho_0
    \psi_\el
    (c, \grad_0 \vect{u}_{\PP}))
    \big),
    \vspace{0.2cm}
    \label{eq:weak_formulation_primal_b}
    \\
    \tens{0}
    & =
    \minus
    \big( \gradL \vect{\xi}, \tens{P}_\PP
    ( c, \gradL \vect{u}_\PP )
    \big)
    +
    \big(\vect{\xi}, \tens{P}_\SEI
    \cdot \normal_\PP
    \big)_{\Gamma_{\text{inter}}}
    \label{eq:weak_formulation_primal_c}
    \\
    \tens{0}
    & =
    \minus
    \big( \gradL \vect{\chi}, \tens{P}_\SEI
    (\gradL \vect{u}_\SEI, \tens{F}_\pl,
    \tens{P}_\Pi (\grad_0 \vect{u}_\SEI,
    \tens{F}_\pl,
    \varepsilon_\pl^\text{eq})
    )
    \big)
    \nonumber
    \\
    &
    \quad
    -
    \big(\vect{\chi}, \tens{P}_\PP
    \cdot \normal_\PP
    \big)_{\Gamma_{\text{inter}}}
    \label{eq:weak_formulation_primal_d}
  \end{empheq}
\end{subequations}
for all test functions~$\varphi \in V$,
$\vect{\xi} \in \vect{V}_\PP^{*}$,
$\vect{\chi} \in \vect{V}_\SEI^{*}$
and~$\tens{P}_\Pi (\grad_0 \vect{u}_\SEI,
\tens{F}_\pl,
\varepsilon_\pl^\text{eq})
=
\Ctensor[\tens{E}_{\el, \logg}]$
with
appropriate function
spaces~$V, \vect{V}_\PP^*$
and~$\vect{V}_\SEI^*$
to guarantee a well-posed formulation~\cite{neff2003comparison}.
The last two spaces include constraints to avoid rigid body motions.
In~\cref{eq:weak_formulation_primal},
$(\something, \something)$
indicates the~$L^2$-inner product for two scalar functions, vectors and
tensors, respectively.
The notation with boundary index~$(\something,
\something)_{\Gamma_{\text{inter}}}$
denotes the interface integral at the interface~${\Gamma_{\text{inter}}}$
regarding~$\normal_\PP$.
\cref{eq:ibvp_e} results in a saddle point problem, which
is condensed
into~\cref{eq:weak_formulation_primal_d}
with the respective projector formulation for the rate-independent
or rate-dependent plastic
approach~\cite{schoof2023efficient}.

\textbf{Discretization.}
Considering an admissible
mesh~$\mathcal{T}_h$
as discretization of the computational domain~$\Omega_\text{h}$,
we apply the isoparametric Lagrangian finite element
method~\cite[Ch.~III~§2]{braess2007finite}
and insert spatial
approximations~$c_h, \mu_h, \vect{u}_{\PP, h}$
and~$\vect{u}_{\SEI, h}$
with the finite dimensional subspaces of the respective function space:
\begin{subequations}
  \begin{alignat}{2}
    &V_h
    &&= {\text{span}} \{ \varphi_i \fdg i=1,\dots, N_\PP \}
    \subset {V},
    \\
    &\vect{V}_{\PP, h}^{*}
    &&= {\text{span}} \{ \vect{\xi}_j \fdg j=1,\dots, 3 N_\PP \}
    \subset \vect{V}_{\PP}^{*},
    \\
    &\vect{V}_{\SEI, h}^{*}
    &&= {\text{span}} \{ \vect{\chi}_k \fdg k=1,\dots, 3 N_\SEI \}
    \subset \vect{V}_{\SEI}^{*}.
  \end{alignat}
\end{subequations}
For the temporal discretization,
we collect all time-dependent coefficients
in the
vector valued function
\begin{align}
  \label{eq:discrete_solution_vector}
  \vect{y} \colon[0, \tfinal] \rightarrow \IR^{(2+3)N_\PP + 3 N_\SEI}, \quad
  t \mapsto \vect{y}(t) =
  \begin{pmatrix}
    \vect{c}_h(t)
    \\
    \vect{\mu}_h(t)
    \\
    \vect{u}_{\PP,h}(t)
    \\
    \vect{u}_{\SEI,h}(t)
  \end{pmatrix}
\end{align}
satisfying
$\tens{M} \ptl_t \vect{y} - \vect{f}(t, \vect{y},
\tens{F}_{\pl, h}, \varepsilon_{\pl, h}^\text{eq})
= \tens{0}$
for
$t \in (0, \tfinal]$
with $\vect{y}(0) = \vect{y}^0$
and the discrete version of the internal
variables~$\tens{F}_{\pl, h},
\varepsilon_{\pl, h}^\text{eq}$~\cite{schoof2023efficient}.
For the temporal discretization of the internal variables,
we use an implicit exponential map~\cite{schoof2023efficient}.
In total,
we get
the space and time discrete problem to go on
from one time~$t_n$
to the next~$t_\npo = t_n + \tau_n$ with $\tau_n> 0$
applying the numerical differential formulation~(NDF)
of linear multistep
methods~\cite{reichelt1997matlab, shampine1997matlab,
shampine1999solving}.
Then,
for given~$\tens{F}_{\pl,h}^n$
and~$\varepsilon_{\pl,h}^{\text{eq},n}$
find the discrete
solution~$\vect{y}^{\npo} \approx \vect{y}(t_{\npo})$
satisfying
\begin{align}
  \label{eq:space_and_time_discretization}
  \alpha_{k_n} \tens{M} \left(\vect{y}^\npo - \vect{\Phi}^n\right)
  - \tau_n
  \vect{f} \left(t_\npo, \vect{y}^\npo,
  \tens{F}_{\pl, h}^{n},
  \varepsilon_{\pl,h}^{\text{eq}, n}
  \right)
  = \tens{0}
\end{align}
with~$\vect{\Phi}^n$,
consisting of solutions on former time
steps~$\vect{y}^n, \dots, \vect{y}^{n-k}$,
and
a constant~\mbox{$\alpha_{k_n}>0$}, which
dependents on the chosen time integration order~$k_n$
at time~$t_n$~\cite[Sect.~2.3]{shampine1997matlab}.
With the computed discrete solution~$\vect{y}^{n+1}$,
we can then update~$\tens{F}_{\pl,h}^{n+1}$
and~$\varepsilon_{\pl,h}^{\text{eq},{n+1}}$
for the new time steps as
explained in~\cite{schoof2023efficient}.
Finally, we solve the nonlinear system in each time step with the
Newton--Raphson method
and follow the adaptive algorithm
in~\cite{castelli2021efficient}.




\section{Numerical Experiments}
\label{sec:results}

\begin{figure}[b]
  \centering
  \includegraphics[width = 0.48\textwidth,
  page=1]{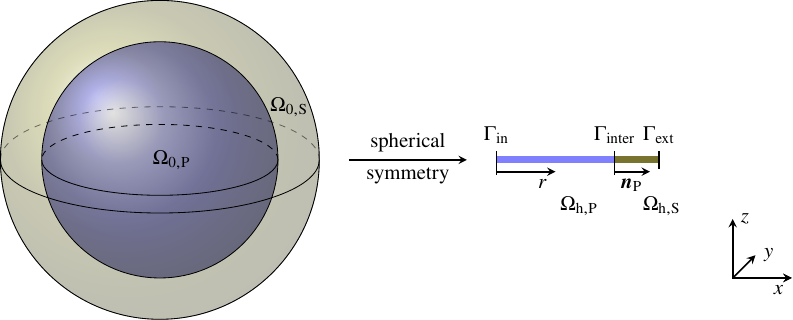}
  \caption{Dimension reduction of a three-dimensional unit sphere with
    surrounded SEI to the one-dimensional interval, based
    on~\cite[Fig.~B.1]{castelli2021numerical}.}
  \label{fig:particle_obstacle_1d}
\end{figure}

\textbf{Simulation Setup.}
We consider a 3D spherical symmetric aSi particle with surrounding
spherical symmetric SEI.
Therefore, we reduce the computational domain due to symmetry assumptions to
the 1D unit interval for the particle~$\Omega_{\text{com}, \PP}$
with additional computational SEI domain~$\Omega_{\text{com}, \SEI}$
as displayed in~\cref{fig:particle_obstacle_1d}.
This means that we have only changes along the radius~$r$,
compare also~\cite[App.~B.2.1]{castelli2021numerical}
and~\cite{kolzenberg2022chemo-mechanical, castelli2021efficient,
schoof2023efficient}.
For the computational setting,
we introduce additional necessary boundary conditions with
no flux and no displacement
\begin{align}
  \vect{N} \cdot \normal_{\PP} = 0,
  \quad
  \vect{u}_\PP = 0
  \quad
  \ton \Gamma_{\text{in}, \tfinal}.
\end{align}
We charge our particle with 1~\si{C}
and discharge with~$\minus 1~\si{C}$
applying the OCV curve stated
in~\cite{kolzenberg2022chemo-mechanical, schoof2023efficient}
with~$c_0 = 0.02$
and a cycling duration of~$0.9~\si{\hour}$.
Further initial conditions are selected as follows:
$\mu_0 = \ptl_c (\rho_0 \psi_\ch(c_0))$,
$u_\PP = r (\lambda_\ch(c_0) - 1)$
and
$u_\SEI = \lambda_\ch(c_0) - 1$
to decrease the number of required Newton steps at the
beginning~\cite{castelli2021efficient}.
All numerical simulations are performed with an isoparametric fourth-order
Lagrangian finite element method
using an integral evaluation through a Gau\ss-Legendre quadrature formula with
six quadrature points in spatial direction and the finite element
library~\texttt{deal.II}~\cite{arndt2023deal}.
Further hardware specifications are given in~\cite{schoof2023efficient}.
Shared memory parallelization is used for assembling of the Newton method.
We solve the linear system with a LU-decomposition.
Due to implementation reasons of the coupled domain,
we use a constant and uniform distributed mesh.
In our experiments, the mesh has around \num{15e3}~degrees of freedom.
The initial time step is chosen as \num{e-8},
the maximal time step as \num{e-3},
$\RelTol_t = \num{e-5}$
and $\AbsTol_t = \num{e-8}$.

\textbf{Numerical Results.}
The parameter setup for the particle can be found in~Table~2
in~\cite{schoof2023efficient}
with further parameters for the SEI taken
from~\cite{kolzenberg2022chemo-mechanical}
with $L_{0, \SEI} = 0.1 L_{0, \PP}$,
$\nu_\SEI = 0.25$,
$E_\SEI = 900\,\si{\mega \pascal}$,
$\sigma_\text{Y} = 49.5\,\si{\mega \pascal}$.
The parameter for viscoplastic plasticity are chosen as
$\dot{\varepsilon}_0 = \num{e-3}\,\si{\per\second}$ or
$\dot{\varepsilon}_0 =
\num{e-4}\,\si{\per\second}$,
$\sigma_{\text{Y}^*} = \sigma_\text{Y}$
and
$\beta = 2.94$.
\begin{figure*}[t]
  \centering
  \includegraphics[width = 0.79\textwidth,
  page=1]{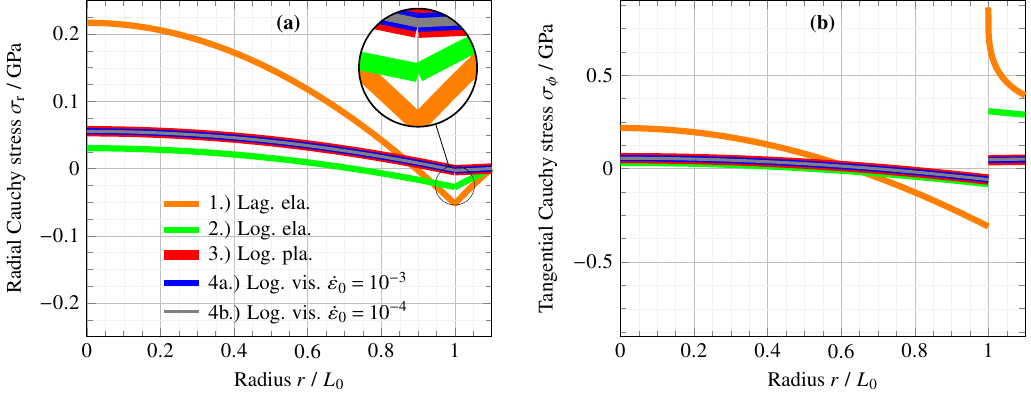}
  \caption{Radial (a) and tangential (b) Cauchy stress over the
  radius ($0 \leq r \leq 1.0$ for $\Omega_{\text{h}, \PP}$, $1.0 \leq r \leq
  1.1$ for $\Omega_{\text{h}, \SEI}$) at the end of the simulation;
  for the
  1.)
  purely elastic (ela.) case with the GSV strain approach (Lag. ela.),
  the
  simulation stopped at $\text{SOC}= 0.34$, see \cref{fig:voltage}(a),
  whereas all
  other cases (2.) elastic, 3.) plastic (pla.) and 4a.), 4b.) viscoplastic
  (vis.))
  for the
  logarithmic (Log.) Hencky strain approach reached the final simulation time.}
  \label{fig:stress}
\end{figure*}
\begin{figure*}[t]
  \centering
  \includegraphics[width = 0.79\textwidth,
  page=1]{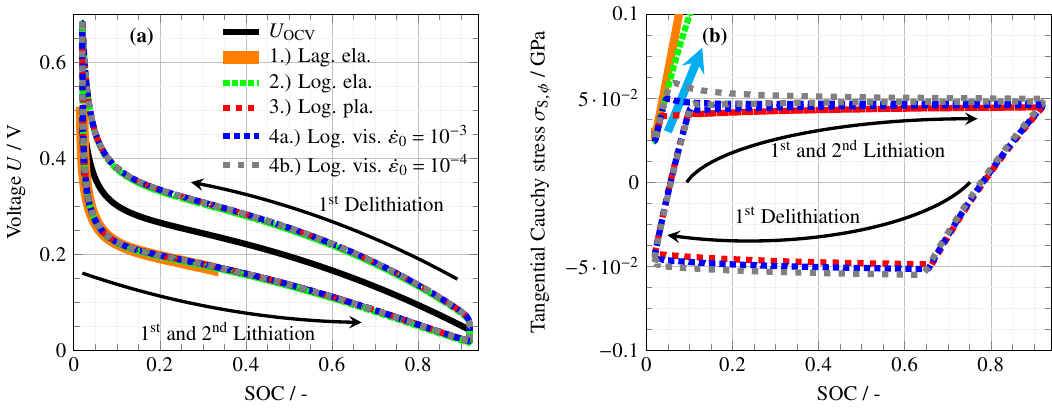}
  \caption{Electrical voltage $U$ over SOC for three half cycles for different
  elastic strain approaches for elastic and plastic cases with aborted GSV
  Lagrangian approach for the elastic case in orange (a)
  and tangential Cauchy SEI
  stress over SOC at
  the particle SEI interface with increasing stress-overrelaxation for
  smaller~$\dot{\varepsilon}_0$
  (blue arrow)
  in the viscoplastic case (b).}
  \label{fig:voltage}
\end{figure*}
We consider in our numerical simulation three half cycles which means two
lithiations and one delithiation.

During the model extension from a single particle setting to the coupled
particle-SEI approach,
we started with a purely elastic SEI with the GSV elastic strain.
However, the numerical simulation stopped around $t = \num{0.32}~\si{\hour}$
and $\text{SOC} = \num{0.34}$, respectively, failing to find an
appropriate Newton update.
\cref{fig:stress}(a)--(b)
and~\cref{fig:voltage}(a)--(b)
shows the numerical results in solid orange.
It is clearly visible in \cref{fig:stress}(b)
that the tangential Cauchy stress of the SEI shows a large
increase at the interface between the particle and the SEI leading to unnatural
behavior and the stop of
the numerical simulation.
On the other hand and also in the purely elastic setting,
the use of the logarithmic strain definition results in a stabilization of the
numerical simulation and suppresses the large increase at the interface~(solid
green).
\cref{fig:stress}(a) shows the required displacement interface condition between
the particle and the SEI.
Due to the application of the exponential map of the plastic time
integration,
the logarithmic strain approach is favorable compared to the GSV strain
approach and is used in all plastic simulations.
\cref{fig:stress}(a)--(b)
and~\cref{fig:voltage}(a)
show almost no change between the rate-independent and rate-dependent
plasticity.
However,
\cref{fig:voltage}(b) shows the typical
stress-overrelaxation for the tangential Cauchy stress of the SEI at the
particle-SEI
interface,
consisting of an
overshooting at the beginning of the first plastification
and followed
by some
relaxation towards the rate-independent plastic results.
A decrease of $\dot{\varepsilon}_0$ from $\num{e-3}\,\si{\per\second}$
to $\num{e-4}\,\si{\per\second}$
results in a higher overshooting~(blue arrow).
An electrical overrelaxation can also be seen in measurements
in~\cite{wycisk2023challenges}
explaining the voltage relaxation after low current
GITT-pulses~\cite{kobbing2024voltage}.
The results fit qualitatively to that ones
in~\cite[Figure~8]{kolzenberg2022chemo-mechanical}.
The purely elastic case has no hysteresis meaning that lithiation and
delithiation result in the same stresses, plotted in dotted green
in~\cref{fig:voltage}(b).
Finally,
\cref{fig:voltage}(a) shows the evaluation of the electric voltage
with the
Butler--Volmer condition in a postprocessing step leading to similar results
compared to the single particle setting~\cite[App.~G]{schoof2023efficient}.




\section{Summary and Conclusion}
\label{sec:conclusion}

In this work,
we consider a thermodynamically consistent chemo-elasto-plastically coupled
model for spherical symmetric aSi particles with surrounded spherical symmetric
SEI using a finite
deformation approach.
We base our theory for the coupling of particle and SEI
on~\cite{kolzenberg2022chemo-mechanical}.
However,
we use the rate-independent and rate-dependent plasticity ansatz
by~\cite{schoof2023efficient}, which is new for a plastic SEI approach.
We apply straightforwardly, for our model extension of the particle-SEI
setup
with plastic SEI approach, the efficient adaptive temporal
algorithm~\cite{castelli2021efficient}
combined with higher order finite element methods on a uniform mesh.
During the development of purely elastic effects,
we discovered that the Green--St-Venant (Lagrangian) strain approach leads to
an abortion during the numerical simulation,
whereas a logarithmic strain approach can overcome this failure
and results in a successful numerical simulation.
The logarithmic approach is also favorable to add plastic deformation due to
the used exponential map of the plastic part of the deformation gradient.
For the rate-dependent plasticity,
the typical stress-overrelaxation can be observed.

In future,
we want to consider
a spatial adaptive algorithm
for the coupled particle-SEI setting
with
a combined discretization and iteration error
estimation~\cite{carraro2018adaptive}
or a residual based error estimator~\cite{schoof2024residual}.


  \vspace{-0.24cm}
  \section*{Declaration of competing interest}
  \vspace{-0.24cm}
  \noindent
  The authors declare that they have no known competing financial interests or
  personal relationships that could have appeared to influence the work in this
  paper.
  \vspace{-0.24cm}
  \section*{CRediT authorship contribution statement}
  \vspace{-0.24cm}
  \noindent
  \textbf{R. Schoof:}
  Methodology, Software, Validation, Formal analysis, Investigation,
  Data Curation, Writing -- original draft, Visualization.
  \textbf{G. F. Castelli:}
  Methodology, Writing -- review \& editing.
  \textbf{W. Dörfler:}
  Conceptualization, Resources, Writing -- review \& editing, Supervision,
  Project administration, Funding acquisition.

  \vspace{-0.24cm}
  \section*{Acknowledgement}
  \vspace{-0.24cm}
  \noindent
  The authors thank L. von Kolzenberg and L. Köbbing
  for intensive and constructive discussions about modeling silicon particles
  as well as A. Dyck, J. Niermann and T. Böhlke for modeling plastic
  effects.
  R.S.
  acknowledges financial
  support by the German Research
  Foundation~(DFG) through the Research Training Group 2218
  SiMET~--~Simulation of Mechano-Electro-Thermal processes in Lithium-ion
  Batteries, project number 281041241.
  \vspace{-0.24cm}

  \section*{ORCID}
  \vspace{-0.24cm}
  \noindent
  R. Schoof:
  \href{https://orcid.org/0000-0001-6848-3844}{0000-0001-6848-3844},
  G. F. Castelli:
  \href{https://orcid.org/0000-0001-5484-6093}{0000-0001-5484-6093},
  W. D\"orfler:
  \href{https://orcid.org/0000-0003-1558-9236}{0000-0003-1558-9236}

  \vspace{-0.24cm}
  \bibliography{./literature/literature}

\end{document}